# Inverse Gram Matrix Methods for Prioritization in Analytic Hierarchy Process: Explainability of Weighted Least Squares Optimization Method

Kevin Kam Fung YUEN*

School of Science, Monash University Malaysia

Jalan Lagoon Selatan, 47500 Bandar Sunway

Selangor Darul Ehsan, Malaysia

Email: kevin.yuen@monash.edu; kevinkf.yuen@gmail.com;

## Abstract

This paper proposes Inverse Gram Matrix (IGM) methods to prioritize Pairwise Reciprocal Matrices (PRMs) in the Analytic Hierarchy Process (AHP). The IGM methods include Pseudo-IGM, Normalized-IGM, and Lagrangian-IGM. Interestingly, the proposed IGM methods achieve the minimum Weighted Least Squares (WLS) error. Since clarity, explainability, usability, and verification of the closed-form solutions of WLS appear to be incomplete in the literature, comprehensive mathematical proofs, detailed computational demonstrations, and intensive simulation verification to extend the prior studies are offered in this study. After a simulation of one million random PRM instances was performed to verify the equivalence of the IGM methods, another simulation of ten thousand random PRM instances was performed to verify that any one of the proposed IGM methods is the exact closed-form solution to the WLS optimization problem. The proposed IGM methods under the WLS framework provide the promising alternatives to Saaty's Eigenvector method in AHP. The open-source R package igm4ahp, along with its user manual, tutorials, and simulation datasets, is publicly available on GitHub and permanently archived on Zenodo.

## 1. Introduction

The first use of pairwise comparisons may be attributed to Ramon Llull in the 13th century [1] [2]. "The Law of Comparative Judgments" was developed for psychological research in [3]. [4] developed pairwise comparisons based on paired ratio scales to establish Analytic Hierarchy Process (AHP) [5] and Analytic Network Process (ANP) [6]. [7, 8] proposed the Cognitive Network Process (CNP) based on the Pairwise Opposite Matrix. The AHP should not be equated with pairwise comparisons [2], as several forms of PCs exist. Despite the increasing use of the AHP, AHP research encompasses several controversial topics, for example, arbitrary hierarchical composition [9-11], conflict with expected utility theory [12, 13], defending against eighteen critics of AHP [14], and rank reversals [15-17].

Prioritization methods for AHP remain among these unsettled, controversial topics. The default and early prioritization method proposed by [4] was the Eigenvector method, and [5] introduced normalization methods including Normalization of the Row Sums, Normalization of Reciprocals of Column Sums, and Arithmetic Mean of Normalized Columns, which were named by [18] as [5] did not explicitly name them. [19] proposed the Direct Least Squares (DLS) and Weighted Least Squares (WLS) methods. [20] proposed the Normalization of Geometric Means, a closed-form solution of the Logarithmic Least Squares (LLS). [21] proposed Enhanced Goal Programming. [22] proposed Fuzzy Programming. [23] proposed Singular Value Decomposition. [24] proposed Maximization of Correlation Coefficients. [25] proposed a linear programming method. [26] proposed the Cosine Maximization Method. [27] proposed the Least Penalty Optimization Prioritization Operator. [28] proposed the Relative Deviation Interconnection method. [29] proposed the Two-Stage Ranking method. Beyond theoretical developments, AHP and MCDM methodologies have seen broad application in real-world decision problems, including land suitability evaluation for organic farming [30] and healthcare decision frameworks such as malignant carcinoma ranking using hesitant VIKOR fuzzy sets [31]. Numerous reviews and comparative studies of prioritization methods exist (e.g., [18, 32-35]), but none of them, including this study, are exhaustive.

This study proposes Inverse Gram Matrix (IGM) methods under the WLS optimization model, providing mathematical clarity and exact closed-form solutions to WLS prioritization. The primary theoretical and computational contribution of this work—deriving exact global WLS minimizers in $O(n^3)$ direct algebraic operations without iterative numerical solvers—is rigorously verified by open-source R code and extensive Monte Carlo simulations. Specifically, a simulation of 1,000,000 random Pairwise Reciprocal Matrices (PRMs) confirms 100% numerical weight equivalence across the three proposed IGM variants (Pseudo-IGM, Normalized-IGM, and Lagrangian-IGM), while a second simulation of 10,000 random PRMs demonstrates that closed-form IGM inversion achieves the exact global WLS minimum over 33,000× faster than non-linear optimization solvers with zero approximation error. The open-source R package (igm4ahp) and simulation datasets are available on GitHub and permanently archived on Zenodo [36]. The package user manual, tutorials, and numerical demonstrations implementing the algorithms and methods in this paper are hosted at https://kkfyuen.github.io/igm4ahp/ .

The rest of this paper is organized as follows. The related work and motivations for this study are presented in Section 2. The Pseudo-IGM, Normalized-IGM, and Lagrangian-IGM methods are illustrated in Section 3. Two numerical examples demonstrate for the usability of IGMs in Section 4. In Section 5, two simulations verify the equivalence of the proposed IGM methods, which produce and demonstrate that they yield the exact and unique closed-form solution to the WLS optimization problem. Conclusions and future research directions are presented in Section 6.

## 2. Related work and Motivation

The rating scores of pairwise comparisons in AHP are encoded in the format of Pairwise Reciprocal Matrix (PRM),

$$A = \{a_{ij}: 0 < a_{ij} = a_{ji}^{-1}, \forall i, \forall j \in (1, \dots, n)\}. \tag{1}$$

$a_{ij}$ is determined by a numeric rating to estimate how many times is object $i$ as important as object $j$. The rating score is usually selected from a 9-point scale. The term "pairwise reciprocal matrix" stems from the reciprocity axiom $0 < a_{ij} = a_{ji}^{-1}$.

To interpret the pairwise matrix, let a set of ideal relative weights (or a priority vector) be $w = [w_1, \dots, w_n]$ such that unit-sum normalization constraint is defined as follows.

$$\sum_{i=1}^{n} w_i = 1. \tag{2}$$

$w$ can be determined by a subjective judgmental PRM denoted by $A = [a_{ij}]$, where a comparison score is

$$a_{ij} \cong \frac{w_i}{w_j}. \tag{3}$$

The ideal PRM $\tilde{A} = \left[\frac{w_i}{w_j}\right]$ relates to the judgment matrix $A = [a_{ij}]$ as follows [18].

$$\tilde{A} = [\tilde{a}_{ij}] = \begin{bmatrix} \frac{w_1}{w_1} & \frac{w_1}{w_2} & \cdots & \frac{w_1}{w_n} \\ \frac{w_2}{w_1} & \frac{w_2}{w_2} & \cdots & \frac{w_2}{w_n} \\ \vdots & \vdots & \ddots & \vdots \\ \frac{w_n}{w_1} & \frac{w_n}{w_2} & \cdots & \frac{w_n}{w_n} \end{bmatrix} \cong \begin{bmatrix} a_{11} & a_{12} & \cdots & a_{1n} \\ a_{21} & a_{22} & \cdots & a_{2n} \\ \vdots & \vdots & \ddots & \vdots \\ a_{n1} & a_{n2} & \cdots & a_{nn} \end{bmatrix} = [a_{ij}] = A, \tag{4}$$

where $a_{ij} = a_{ji}^{-1}$, and $a_{ii} = 1$, $i, j = 1, \dots, n$. Here, $A$ is the PRM derived from expert judgment. A prioritization method derives $w$ from $A$. Examples of prioritization methods are outlined in the introduction section. This research focuses on development of IGM methods closely related to the Weighted Least Squares (WLS) method. Chu et al. [19] proposed the Direct Least Squares (DLS) and WLS methods. The DLS method minimizes the total variance between $a_{ij}$ and $\tilde{a}_{ij} = \frac{w_i}{w_j}$.

$$\begin{aligned} &Min \quad \sum_{i=1}^{n}\sum_{j=1}^{n}\left(a_{ij} - \frac{w_i}{w_j}\right)^2 \\ &\text{Subject to} \sum_{i=1}^{n} w_i = 1,\ w_i > 0, i = 1,2,\dots,n\,. \end{aligned} \tag{5}$$

Chu et al. [19] indicated that the above non-linear optimization problem has no simple tractable form, may have multiple solutions, and is very difficult to solve numerically. Bozóki [37] showed that the optimization problem can be reduced to solving a system of polynomial equations and concluded that the DLS solution is not unique in general, with the number of solutions potentially equal to or up to twice as the matrix dimension. Carrizosa and Messine [38] indicated that such optimization problems are occasionally multimodal and proposed an exact global optimization algorithm to obtain globally optimal weights. Chu et al. [19] modified the objective function and proposed the Weighted Least Squares (WLS) method in the form:

$$Min \quad \sum_{i=1}^{n}\sum_{j=1}^{n}(w_i - a_{ij}w_j)^2 \tag{6}$$

$$\text{Subject to} \sum_{i=1}^{n} w_i = 1,\ w_i > 0, i = 1,2,\dots,n \tag{7}$$

Regarding the closed-form solution, Blankmeyer [39] indicated that the objective function may be expressed in quadratic form $v^T C^{-1} v$, where a square matrix $C$ has the diagonal elements as below.

$$c_{ii} = (n-1) + \sum\nolimits_{j, j \neq i} a_{ji}^2, \qquad i = 1, \cdots, n, \tag{8}$$

and the off-diagonal elements:

$$c_{ij} = -a_{ij} - a_{ji}, \qquad i, j = 1, \cdots, n, \tag{9}$$

The solution of the closed form was expressed as.

$$v = \frac{C^{-1}e}{e^T C^{-1}e} \tag{10}$$

As shown in Eq. (10) [19], if $e$ is placed after a matrix, $e$ denotes a column vector of ones, *e.g., $e = (1, \cdots, 1)^T$*. in this paper, e.g. Eq. (23), if $e$ is placed before a matrix, $e$ denotes a row vector of ones of the size with a length equal to the column dimension of the matrix. Eqs.(8)-(9) were first mentioned by Chu et al. [19], who did not provide clear and correct details explaining how they came from. To bridge the research gap, the derivations of Eqs.(8)-(9) are further explained and illustrated in Sections 3 and 4. Eq.(10) was proposed by Blankmeye [39], who [39] highlighted a major limitation: Eqs.(8)-(10) fail for a perfectly consistent PRM because $C$ is notpositive semi-definite and therefore non-invertible. To bridge the research gap, the proposed Inverse Gram Matrix (IGM) methods presented in the next section resolve this non-invertible problem. In addition, Chu et al. [19] and Blankmeyer [39] did not present relevant proofs of the closed-form solution, step-by-step computational demonstration, or simulation verifications, all of which are addressed in this study. Consequently, prior literature reviews and and comparative studies (e.g., [18, 32-35]), merely acknowledged the contributions from Chu et al. [19] and Blankmeyer [39] without utilizing explicit closed-form solutions that bypass numerical optimization solvers. The objective of this study is to establish and promote the proposed IGM methods for exact WLS prioritization.

## 3. Inverse Gram Matrix Methods

### 3.1. Pseudo Inverse Gram Matrix Method

A design matrix $D$ is a coefficient matrix corresponding to $w^T$. $D$ can be constructed from $A$ shown in Eq. (4) via Algorithm 1. To transform the ratio equality $a_{ij} \cong \frac{w_i}{w_j}$ into linear form, the weighted differential equation is expressed as.

$$w_i - a_{ij}w_j = 0, i \neq j, \forall i, j \in \{1, \dots, n\}. \tag{11}$$

Based on the element forms in Eqs. (2) and (11), linear system is expressed in matrix form as follows.

$$Dw^T = q^T, \tag{12}$$

$$q^T = \begin{bmatrix} \mathbf{0}^T \\ 1 \end{bmatrix}. \tag{13}$$

By Algorithm 1, $A$ is converted into $D$ as follows.

$$D = \begin{bmatrix} \bar{D} \\ e \end{bmatrix} = [d_{kj}] = \begin{bmatrix} 1 & -a_{12} & 0 & \cdots & 0 & 0 \\ -a_{21} & 1 & 0 & \cdots & 0 & 0 \\ 1 & 0 & -a_{13} & \cdots & 0 & 0 \\ -a_{31} & 0 & 1 & \cdots & 0 & 0 \\ \vdots & \vdots & \vdots & \vdots & \vdots & \vdots \\ 1 & 0 & 0 & \cdots & 0 & -a_{1n} \\ -a_{n1} & 0 & 0 & \cdots & 0 & 1 \\ \vdots & \vdots & \vdots & \vdots & \vdots & \vdots \\ 0 & 0 & 0 & \cdots & 1 & -a_{n-1,n} \\ 0 & 0 & 0 & \cdots & -a_{n,n-1} & 1 \\ 1 & 1 & 1 & \cdots & 1 & 1 \end{bmatrix}. \tag{14}$$

**Algorithm 1 ($D = \Pi(A)$: Converting A into D).**

**Input**: A Pairwise Reciprocal Matrix $A = [a_{ij}]$

**Step 1:** Initialize the design matrix without weight sum unity function: $\bar{D} = \mathbf{0}_{n(n-1),n}$

**Step 2:** initialize the first-row index for *D*: $k = 1$

**Step 3:** Fill elements into *D*

For *i* from 1 to *n* – 1*:*

For *j* from (*i*+*1*) to *n:*

$\bar{D}_{k,i} = 1$

$\bar{D}_{k,j} = -a_{ij}$

$k = k + 1$ # move the index to the next row of *D*

$\bar{D}_{k,i} = -a_{ji}$

$D_{k,j} = 1$

$k = k + 1$ # move the index to the next row of *D*

**Step 4:** Insert *e*, a vector of 1, into the last row of *D*: $D = \begin{bmatrix} \bar{D} \\ e \end{bmatrix}$

**Return:** $D$.

An example of *D* derived from a $4 \times 4$ PRM is shown in Eq.(36). The augmented design matrix *D* consists of two blocks: the pairwise differential block $\bar{D}$ and the normalization row vector $e$. The objective differential block $\bar{\mathbf{D}} \in \mathbb{R}^{n(n-1)\times n}$ comprises the first $n(n-1)$ rows of $D$ and corresponds to the pairwise coefficients in Eq.(11). The normalization row vector of ones, *e,* forming the bottom row of *D* enforces the unit-sum constraint in Eq. (2).

$q^T$ is a column vector (the transpose of row vector *q*) with its last element equal to 1 and all preceding elements equal to 0, where $\mathbf{0}^T = [0, \dots, 0]^T$ has a length matching the number of pairwise differential equations. For Eq.(11), the row dimension of $Dw^T$ is $n(n-1) + 1$.

The sequence order of equations in *D* is not important as there are many forms of row-equivalent matrices, i.e., a new matrix swapping rows of a matrix is row equivalent to the original one. Algorithm 1 provides one of examples to generalize the coefficient matrix, or design matrix of a system of linear equations. The gram matrix *G* of the form below is the transpose of *D* (i.e., $D^T$) multiples $D$.

The order of equations in *D* does not affect the solution, as row-swapped matrices remain row-equivalent. Algorithm 1 provides a systematic procedure to construct the coefficient matrix, or the design matrix for the linear system. The Gram matrix *G* is formed by premultiplying *D* by its transpose $D^T$:

$$G = [g_{ij}] = D^T D. \tag{15}$$

*G* is symmetric. To guarantee that the Gram matrix $G = D^T D$ is non-singular and strictly invertible, the rank of the design matrix *D* must be established. The following lemma proves that for any valid Pairwise Reciprocal Matrix (PRM), the design matrix *D* constructed by Algorithm 1 possesses full column rank.

**Lemma 1 (Full Column Rank of Design Matrix $D$).**

For any Pairwise Reciprocal Matrix $A = [a_{ij}]_{n\times n}$ with strictly positive entries $a_{ij} > 0$, the design matrix $D \in \mathbb{R}^{(n(n-1)+1)\times n}$ constructed by Algorithm 1 has full column rank, *i.e.*, $\text{rank}(D) = n$. Consequently, the Gram matrix $G = D^T D$ is symmetric positive-definite and strictly invertible ($\det(G) > 0$).

**Proof.**

Partition the design matrix $D = \begin{bmatrix} \overline{D} \\ e \end{bmatrix}$ into its objective differential block $\overline{D} \in \mathbb{R}^{n(n-1)\times n}$ and its normalization row vector $e = [1,1,\dots,1] \in \mathbb{R}^{1\times n}$.

To establish $\text{rank}(D) = n$, consider an arbitrary vector $u^T = [u_1, u_2, \dots, u_n]^T \in \mathbb{R}^n$ in the kernel (null space) of $D$, satisfying the homogeneous condition $Du^T = \mathbf{0}$:

$$\text{null}(D) = \{u^T \in \mathbb{R}^n : Du^T = \mathbf{0}\} = \{\mathbf{0}\} \tag{16}$$

(Note: The zero vector $\mathbf{0}$ is used here strictly to define the matrix kernel for column rank evaluation, which is distinct from the inhomogeneous right-hand side $q^T = [0,\dots,0,1]^T$ in Eq. (12) used to solve for priority weights).

Expanding $Du^T = \mathbf{0}$ block-wise yields:

$$\begin{bmatrix} \overline{D} \\ e \end{bmatrix} u^T = \begin{bmatrix} \mathbf{0} \\ 0 \end{bmatrix} \Leftrightarrow \begin{cases} \overline{D}u^T = \mathbf{0} \\ eu^T = 0 \end{cases} \tag{17}$$

1. Pairwise Differential Block ($\overline{D}u^T = 0$): For every distinct pair $i \neq j$,

$$\begin{aligned} u_i - a_{ij}u_j &= 0 \\ u_i &= a_{ij}u_j, \quad \forall i \neq j \end{aligned} \tag{18}$$

2. Kernel Constraint ($eu^T = 0$): The bottom row enforces:

$$\sum_{k=1}^{n} u_k = 0 \tag{19}$$

Expressing all components $u_i$ in terms of $u_1$ via Eq. (18) gives $u_i = a_{i1}u_1$ for $i = 1, \dots, n$ . Substituting these into Eq. (19) yields:

$$\begin{aligned} \sum_{i=1}^{n} u_i &= u_1 + u_2 + u_3 + \cdots + u_n = 0 \\ \sum_{i=1}^{n} u_i &= (a_{11}u_1) + (a_{21}u_1) + (a_{31}u_1) + \cdots + (a_{n1}u_1) = 0 \\ \sum_{i=1}^{n} u_i &= u_1 + a_{21}u_1 + a_{31}u_1 + \cdots + a_{n1}u_1 = 0 \\ \sum_{i=1}^{n} u_i &= u_1(1 + a_{21} + a_{31} + \cdots + a_{n1}) = 0 \\ \sum_{i=1}^{n} u_i &= u_1\left(1 + \sum_{i=2}^{n} a_{i1}\right) = 0 \end{aligned} \tag{20}$$

Since $a_{i1} > 0$ for all $i$, the multiplier satisfies

$$\left(1 + \sum_{i=2}^{n} a_{i1}\right) > 1 > 0. \tag{21}$$

Therefore, the equation requires $u_1 = 0$, which immediately forces:

$$u_i = a_{i1}(0) = 0, \quad \forall i = 1, \dots, n \tag{22}$$

Thus, $Du^T = \mathbf{0} \Longrightarrow u^T = \mathbf{0}$, establishing that $\text{null}(D) = \{0\}$ and $\text{rank}(D) = n$ (full column rank) for any positive PRM $A$. Because $D$ has full column rank, $G = D^T D$ is symmetric positive-definite, non-singular, and strictly invertible ($\det(G) > 0$). QED.

Having established in Lemma 1 that the Gram matrix $G = D^T D$ is strictly positive-definite and non-singular ($\det(G) > 0$), the closed-form priority weight derivation for the Pseudo Inverse Gram Matrix (PIGM) method is formally presented in Theorem 1.

**Theorem 1 (**Pseudo-IGM**).**
The matrix-form solution of Pseudo-IGM is given by:

$$w = \frac{eG^{-1}}{eG^{-1}e^T} = \frac{e(D^T D)^{-1}}{e(D^T D)^{-1}e^T} \tag{23}$$

As $G$ and $G^{-1}$ are symmetric, the column sums $eG^{-1}$ (a row vector) and row sums $e^T G^{-1}$ (a column vector) yield identical numerical values. The scalar denominator $eG^{-1}e^T$ represents the sum of all elements in $eG^{-1}$ (the grand sum of $G^{-1}$), ensuring that the entries of $w$ sum to unity (unit-sum normalization).

**Proof:**

From Eq. (12), we have

$$Dw^T = q^T. \tag{24}$$

Premultiplying both sides by $D^T$ yields

$$D^T D w^T = D^T q^T. \tag{25}$$

Since

$$D^T q^T = e^T, \tag{26}$$

substituting this into to Eq. (25) yields the unnormalized linear system.

$$D^T D v^T = e^T, \tag{27}$$

where $v$ denotes the unnormalized priority weight vector. As equality of Eq. (2) may not be preserved, the normalized function is as follows.

$$w^T = \frac{v^T}{ev^T}. \tag{28}$$

Thus,

$$v^T = (D^T D)^{-1} e^T = G^{-1} e^T, \text{ or} \tag{29}$$

$$v = e(D^T D)^{-1} = eG^{-1}, \tag{30}$$

Normalizing $v$ yields:

$$w = \frac{v}{ve^T} = \frac{eG^{-1}}{eG^{-1}e^T} = \frac{e(D^TD)^{-1}}{e(D^TD)^{-1}e^T}. \tag{31}$$

Thus Eq. (23) is hold. QED.

The unnormalized matrix-form solution of Pseudo-IGM is obtained by applying the standard Moore-Penrose pseudoinverse $D^+ = (D^TD)^{-1}D^T = G^{-1}D^T$ to the linear system $D\boldsymbol{w}^T = q^T$:

$$\mathbf{v}^T = D^+q^T = (D^TD)^{-1}D^Tq^T = (D^TD)^{-1}e^T = G^{-1}e^T \tag{32}$$

As the design matrix $D$ has full column rank ($\mathrm{rank}(D) = n$, Lemma 1), the Gram matrix $G = D^TD$ is strictly invertible ($\det(G) > 0$). Consequently, the expression $G^{-1}D^T$ is mathematically identical to the standard Moore-Penrose pseudoinverse $(D^TD)^{-1}D^T$. Furthermore, because $D^Tq^T = e^T$(Eq.(26)) (since $q^T = [0, \ldots, 0, 1]^T$extracts the bottom normalization row $e = [1, 1, \ldots, 1]$), premultiplying by $D^+$ reduces algebraically to $G^{-1}e^T$. As $G^{-1}$is symmetric, $\boldsymbol{v} = eG^{-1}$represents the row sums of $G^{-1}$. The calculation steps of Pseudo-IGM are summarized as follows. Firstly, a design matrix $D$ is converted from a PRM $A$ by Algorithm 1. Secondly, a gram matrix is computed by Eq. (15). Thirdly, inverse gram matrix is computed by a general inverse method. Finally, $w$ is computed by Eq. (23).

## 3.2. Normalized Inverse Gram Matrix Method

To avoid constructing the explicit design matrix $D$ via Algorithm 1, an alternative approach is to directly convert a PRM $A = [a_{ij}]$ into a Gram matrix in element form, as presented in Theorem 2.

**Theorem 2 (Element form of Gram matrix).**

The elements of the Gram matrix $G = [g_{ij}]_{n\times n} = D^TD$ derived from an $n \times n$ Pairwise Reciprocal Matrix $A = [a_{ij}]_{n\times n}$ are given by

$$g_{ij} = \begin{cases} (n-1) + \sum_k a_{kj}^2 & i = j \\ 1 - a_{ij} - a_{ji} & i \neq j \end{cases} \tag{33}$$

or equivalently,

$$g_{ij} = \begin{cases} n + \sum_{k\neq j} a_{kj}^2 & i = j \\ 1 - a_{ij} - a_{ji} & i \neq j \end{cases} \tag{34}$$

**Proof.**

To illustrate the matrix structure intuitively, consider a $4 \times 4$ PRM as an initial concrete case

$$A = [a_{ij}] = \begin{bmatrix} 1 & a_{12} & a_{13} & a_{14} \\ a_{21} & 1 & a_{23} & a_{24} \\ a_{31} & a_{32} & 1 & a_{34} \\ a_{41} & a_{42} & a_{43} & 1 \end{bmatrix} \tag{35}$$

Based on Algorithm 1, the design matrix $D \in \mathbb{R}^{13\times 4}$ is:

$$D = \begin{bmatrix} 1 & -a_{12} & 0 & 0 \\ -a_{21} & 1 & 0 & 0 \\ 1 & 0 & -a_{13} & 0 \\ -a_{31} & 0 & 1 & 0 \\ 1 & 0 & 0 & -a_{14} \\ -a_{41} & 0 & 0 & 1 \\ 0 & 1 & -a_{23} & 0 \\ 0 & -a_{32} & 1 & 0 \\ 0 & 1 & 0 & -a_{24} \\ 0 & -a_{42} & 0 & 1 \\ 0 & 0 & 1 & -a_{34} \\ 0 & 0 & -a_{43} & 1 \\ 1 & 1 & 1 & 1 \end{bmatrix} \tag{36}$$

For $G = [g_{ij}] = D^T D$, a diagonal element, $g_{ii}$ or $g_{jj}$, is the sum of the squared elements of the column $j$ in $D$, (or row $i$ in $D^T$). When $j$ is 1, the first diagonal element is

$$\begin{aligned} g_{11} &= 3(1)^2 + (-a_{21})^2 + (-a_{31})^2 + (-a_{41})^2 + 1^2 \\ &= 3(1) + (1 + a_{21}^2 + a_{31}^2 + a_{41}^2) \\ &= (4-1) + \sum\nolimits_{k=1}^{4} (-a_{k1})^2 \\ &= 4 - 1 + \sum_{k=1}^{4} a_{k1}^2 \end{aligned} \tag{37}$$

More generally, for an arbitrary $n \times n$ PRM, column $j$ of $D \in \mathbb{R}^{(n(n-1)+1)\times n}$ contains $(n-1)$ entries equal to $+1$ (from rows $w_j - a_{jk}w_k = 0$), $(n-1)$ entries equal to $-a_{kj}$ for $k \neq j$ (from rows $w_k - a_{kj}w_j = 0$), and 1 entry equal to $+1$ in the bottom normalization row $e$. Summing the squared elements of column $j$ yields the general diagonal entry for any dimension $n \geq 3$.

$$\begin{aligned} g_{ii} = g_{jj} &= \sum\nolimits_{k=1}^{n(n-1)+1} D_{kj}^2 \\ &= (n-1)(1)^2 + \sum\nolimits_{k\neq j} \left(-a_{kj}\right)^2 + 1^2 \\ &= (n-1) + \sum\nolimits_{k} a_{kj}^2 = (n-1) + \left(a_{jj}^2 + \sum\nolimits_{k\neq j} a_{kj}^2\right) \\ &= n - 1 + 1 + \sum\nolimits_{k\neq j} a_{kj}^2 = n + \sum\nolimits_{k\neq j} a_{kj}^2, \qquad j = 1, \cdots, n. \end{aligned} \tag{38}$$

For the off-diagonal entries, when $i = 1$ and $j = 2$,

$$g_{12} = (1)(-a_{12}) + (-a_{21})(1) + (1)(1) = 1 - a_{12} - a_{21}. \tag{39}$$

More generally, for any two distinct columns $i \neq j$ in $D$, only three rows yield non-zero column products: row $w_i - a_{ij}w_j = 0$ contributes $(+1)\left(-a_{ij}\right) = -a_{ij}$; row $w_j - a_{ji}w_i = 0$ contributes $\left(-a_{ji}\right)(+1) = -a_{ji}$; and the bottom row $e$ contributes $(1)(1) = 1$. All other $n(n-1) - 2$ rows contain zero in column $i$ or column $j$. Summing column inner products yields the general off-diagonal entry for any dimension $n \geq 3$:

$$g_{ij} = \sum_{k=1}^{n(n-1)+1} D_{ki}D_{kj} = (-a_{ij}) + (-a_{ji}) + 1 = 1 - a_{ij} - a_{ji} \tag{40}$$

Therefore, the element forms of $G = D^T D$, i.e., Eqs.(33) and (34), hold for any general dimension $n \geq 3$.

QED.

Note: Eqs.(33) and (34) are mathematically equivalent representations of the diagonal entries $g_{jj}$. In Eq. (33) , the summation covers all indices $k = 1, \ldots, n$ (including $a_{jj}^2 = 1^2 = 1$), whereas in Eq. (34), $a_{jj}^2 = 1$ is extracted to simplify $n - 1 + 1 = n$ with index restriction $k \neq j$. In both equations, $n$ denotes the matrix dimension of $A \in \mathbb{R}^{n \times n}$. The only difference between Eq.(34) and the matrix entries $c_{ij}$ in Eqs.(8)-(9) is the addition of "plus one", i.e.,

$$g_{ij} = c_{ij} + 1 \tag{41}$$

This difference arises because the unit-sum normalization constraint ($\sum_{i=1}^{n} w_i = 1$) was omitted from the coefficient matrix in earlier formulations (e.g., Chu et al.[19] and Blankmeyer [39] ). Without this normalization row $e$, the design matrix $D$ reduces to the differential block $\bar{D}$, and the right-hand side vector $q^T$ becomes the zero vector $0^T$. Finding a unique, non-trivial priority vector $w \neq 0$ from the resulting homogeneous system $\bar{D}w^T = 0^T$ leads to $Cw^T = 0^T$. For a perfectly consistent PRM, $\det C = 0$, rendering $C$ non-invertible and preventing a direct solution via $C^{-1}$. To resolve this non-invertibility problem, the Normalized Inverse Gram Matrix (NIGM) method is presented in Theorem 3.

**Theorem 3** (NIGM).

The matrix-form solution of Normalized-IGM (NIGM) is given by:

$$w = \frac{e\tilde{v}}{e\tilde{v}e^T} = \frac{e\tilde{G}^{-1}}{e\tilde{G}^{-1}e^T} = \frac{e(\bar{G}+r)^{-1}}{e(\bar{G}+r)^{-1}e^T} \quad , \tilde{g}_{ij} \in \tilde{G} \,, \tag{42}$$

where $E = e^T e \in \mathbb{R}^{n \times n}$ is the $n \times n$ matrix of all ones, $r \in \mathbb{R}$ ($r \neq 0$), and the elements $\tilde{g}_{ij}$ of $\tilde{G}_r = \bar{G} + rE$ are defined as:

$$\tilde{g}_{ij} = r + \bar{g}_{ij} = \begin{cases} r + (n-1) + \sum_{k \neq i} a_{kj}^2 & i = j \\ r - a_{ij} - a_{ji} & i \neq j \end{cases} , r \in \mathbb{R}, r \neq 0, \forall i, \forall j. \tag{43}$$

Note: The Eq.(43) can be transformed as $\bar{g}_{ij} = \tilde{g}_{ij} - r$.

**Proof.**

Using Theorem 1, Eq.(30) is expressed as.

$$v = e(D^T D)^{-1} = e\left(\begin{bmatrix}\bar{D}\\ e\end{bmatrix}^T \begin{bmatrix}\bar{D}\\ e\end{bmatrix}\right)^{-1} = e(\bar{D}^T\bar{D} + e^T e)^{-1} = eG^{-1} \tag{44}$$

If $r = 0$, the normalization row $e$ is omitted from $D$, reducing $D$ to $\bar{D}$. Thus,

$$\bar{v} = e(\bar{D}^T\bar{D})^{-1} = e\bar{G}^{-1}. \tag{45}$$

Similar to the derivation in Theorem 2, the elements $\bar{g}_{ij}$ of $\bar{G} = \bar{D}^T\bar{D}$ are

$$\bar{g}_{ij} = \begin{cases} (n-1) + \sum_{k \neq i} a_{kj}^2 & i = j \\ -a_{ij} - a_{ji} & i \neq j \end{cases}, \forall i, \forall j. \tag{46}$$

The pseudo-inverse method in Eq. (45) provides the underlying derivation for Blankmeyer's formulation [39] in Eq. (10). $\bar{G}$ in this study is equivalent to $C$ in [39] ,i.e., $\bar{G} = C$. The variant notations of $\bar{G}$, $\tilde{G}$ or $G$ are used in this study to highlight the term "Gram". However, when a PRM is perfectly consistent, $\bar{G}$ (or $C$) is singular ($\det C = 0$) and not invertible, so setting $r = 0$ is not viable. To guarantee invertibility for any PRM, a scalar parameter $r \neq 0$ is introduced under two parameter conditions:

For $r > 0$, constructing $D$ with a scaled normalization row $\sqrt{r}e$ yields:

$$\tilde{v} = e(\tilde{D}^T\tilde{D})^{-1} = e\left(\begin{bmatrix} \bar{D} \\ e\sqrt{r} \end{bmatrix}^T \begin{bmatrix} \bar{D} \\ e\sqrt{r} \end{bmatrix}\right)^{-1} = e\tilde{G}^{-1}. \tag{47}$$

For $r < 0$,

$$\tilde{v} = e(\tilde{D}^T\tilde{D})^{-1} = e\left(\begin{bmatrix} \bar{D} \\ e\sqrt{|r|} \end{bmatrix}^T \begin{bmatrix} \bar{D} \\ -e\sqrt{|r|} \end{bmatrix}\right)^{-1} = e\tilde{G}^{-1}. \tag{48}$$

Following the element expansion steps in Theorem 2, Eqs.(42) and (43) hold. QED.

Eqs. (44)-(45) and the proof steps in Theorem 2 show how Eqs.(8)-(10) are derived, which are not mentioned in [19] and [39]. If $A$ is a perfectly consistent PRM, the remedy for non-invertibility problem using Eqs.(8)-(10) in [39] is shown as below.

$$\tilde{G} = \bar{G} + rE = C + rE = [\bar{g}_{ij}] + rE, r \neq 0 \tag{49}$$

To resolve the non-invertibility problem in Eqs. (8)-(10) identified by Blankmeyer [39] for perfectly consistent Pairwise Reciprocal Matrices (PRMs), a non-zero constant $r$ is simply added to $G$ or $C$. The formal theoretical properties of this parameter and its role in motivating the Lagrangian formulation are established in Remark 1.

**Remark 1 (Invariance under Parameter $r$ and Motivation for Lagrangian-IGM).**

It is important to note that the constant $r \in \mathbb{R}$ ($r \neq 0$) serves exclusively to guarantee the algebraic invertibility of the Gram matrix $G_r = G + r\mathbf{E}$ (where $\mathbf{E} = e^T e \in \mathbb{R}^{n \times n}$ is the $n \times n$ matrix of all ones) when a PRM is perfectly consistent. Because the priority weight vector$w$ is obtained via post-hoc vector normalization:

$$w = \frac{v}{ve^T} = \frac{e(G + r\mathbf{E})^{-1}}{e(G + r\mathbf{E})^{-1}e^T} \tag{50}$$

the scalar constant $r$ cancels out completely during division. Consequently, any non-zero real scalar $r$ produces the exact same normalized weight vector $w$. The parameter $r$ does not alter the priority weights, nor does it affect the Consistency Index (CI), Consistency Ratio (CR), rank ordering, or ratio-scale properties of the decision model. Setting $r = 1$(or $r = 0$ for inconsistent PRMs where $G$ is already non-singular) is recommended for computational simplicity.

Because the parameter $r$ in Normalized-IGM (N-IGM) cancels out during normalization, adjusting $r$ cannot eliminate the requirement for post-hoc vector normalization. This raises a fundamental theoretical question: *Is it possible to formulate a matrix inversion model where explicit vector normalization is not required?* The next two subsections address this issue directly by establishing the Lagrangian-IGM method (Section 3.3) and proving its global optimality and uniqueness (Section 3.4).

## 3.3. Lagrangian Inverse Gram Matrix Method

Chu et al. [19] presented Lagrange expression for Eq.(6) as follows. (Note: The notations of Eqs. (51) - (58) are directly quoted from [19] without modification for the purpose of comparison to this study. They are treated as local variables in [19].) The set of notations does not appear to be mentioned or used in the literature, although some parts are correct and have merit for further investigation and use.

$$S' = \left(\sum_{i=1}^{n}\sum_{j=1}^{n}\left(w_i - a_{ij}w_j\right)^2\right) + 2\lambda\sum_{i=1}^{n} w_i \tag{51}$$

Differentiating Eq.(51) with respect to $w_m$ yields the following expression.

$$\sum_{i=1}^{n}(a_{im}w_m - w_i)\,a_{im} - \sum_{j=1}^{n}\left(a_{mj} - w_m\right) + \lambda = 0, \qquad m = 1, \dots, n. \tag{52}$$

Based on Eq.(52), Chu et al. [19] presented the matrix form as follow.

$$Bw = m, \tag{53}$$

where $B$ is a $(n+1)\times(n+1)$ matrix with elements $b_{ij}$ defined as follows.

$$b_{ii} = \bar{g}_{ij} = c_{ii} = (n-1) + \sum_{j\neq i} a_{ji}^2, \qquad i, j = 1, \cdots, n, \tag{54}$$

$$b_{ij} = \bar{g}_{ij} = c_{ij} = -a_{ij} - a_{ji}, \qquad i, j = 1, \cdots, n, \tag{55}$$

$$b_{k,n+1} = b_{n+1,k} = 1, \qquad k = 1, \cdots, n, \tag{56}$$

$$m = [0, \dots, 0, 1]^T, \tag{57}$$

$$w = [w_1, \dots, w_n, \lambda]^T. \tag{58}$$

Chu et al. [19] did not provide a rigorous proof or detailed procedure for solving the linear system in Eq.(53). Furthermore, the derivation connecting the Lagrangian function in Eq. (51) to its partial derivatives in Eq. (52), as well as the transition from the scalar equations in Eq.(52) to the matrix representation in Eq.(53), remains unclarified in their original work. Although this study establishes that solving Eq.(53) yields the exact closed-form solution for Weighted Least Squares (WLS), Chu et al. [19] did not propose a solution methodology for this system and introduced several algebraic oversights in their intermediate steps. Specifically, the constraint term in Eq.(51) does not conform to the standard definition of the Lagrange multiplier theorem, as the essential unitary constant $(-1)$ is omitted; the mathematically corrected formulation is established in Eq.(59) of this study. Although Chu et al. [19] omitted the partial derivative of Eq.(51) with respect to λ, evaluating this derivative yields $\sum_{i=1}^{n} w_i = 0$, , which erroneously forces the sum of weights to zero rather than enforcing the required normalization constraint $\sum_{i=1}^{n} w_i = 1$. Additionally, the scalar expression in Eq.(52) does not directly map to the matrix form in Eq.(53); for instance, squared entries $a_{im}^2$ persist in the off-diagonal elements of Eq.(52) , whereas the off-diagonal elements $b_{ij}$ defined in Eq. (55) omit these $a_{im}^2$ terms without algebraic justification. To resolve this inconsistency, the mathematically sound partial differentiation is presented in Eq.(60). To provide a rigorous Lagrangian formulation that guarantees the minimum WLS objective error, the Lagrangian-IGM method is formally established in Theorem 4 based on the theorems derived herein.A comprehensive side-by-side comparison between the original formulation of Chu et al. (Eq. (51)) and the corrected Lagrangian formulation (Eq. (59)) is summarized in Table 1.

Table 1. Comparison between Chu et al. [19] and Corrected Lagrangian-IGM

| Property | Chu et al. [19] | Corrected Lagrangian-IGM | Theoretical Resolution |
|---|---|---|---|
| Lagrangian Formulation | $S' = \sum_{i,j}(w_i - a_{ij}w_j)^2 + 2\lambda\sum_i w_i$ (Eq. (51)) | $y = \sum_{i,j}(w_i - a_{ij}w_j)^2 + \bar{\lambda}\left(\sum_i w_i - 1\right)$ (Eq. (59)) | Restores the essential $-1$ constraint constant. |
| Derivative w.r.t. $\lambda$ | $\frac{\partial S'}{\partial \lambda} = 2\sum_{i=1}^{n} w_i = 0 \Longrightarrow \sum_{i=1}^{n} w_i = 0$ | $\frac{\partial y}{\partial \lambda} = \sum_{i=1}^{n} w_i - 1 = 0 \Longrightarrow \sum_{i=1}^{n} w_i = 1$ (Eq. (61)) | Correctly enforces normalization instead of zero. |
| Off-Diagonal Derivative Terms | Retains $a_{im}^2$ terms (Eq. (52)) | Eliminates $a_{im}^2$ via $a_{ii} = 1$ identity (Eqs. (60),(65)-(70) ) | Resolves algebraic mismatch with off-diagonal $g_{ij} = -a_{ij} - a_{ji}$. |
| Matrix Representation | Matrix conversion unproven (Eq. (53)) | Exact map to $\overline{G}[w,\lambda]^T = q^T$ (Eq. (62)) | Guarantees non-singular closed-form solution via $[w^*, \lambda^*] = q\overline{G}^{-1}$ (Eq. (64)). |

To establish the direct theoretical connection between the Weighted Least Squares (WLS) optimization problem and the Inverse Gram Matrix framework, Theorem 4 derives the Karush-Kuhn-Tucker (KKT) first-order stationary conditions of WLS. This demonstrates that the closed-form solution of Lagrangian-IGM (LIGM) is the exact analytical representation of WLS, forming the theoretical foundation for all three equivalent IGM variants (PIGM, NIGM, and LIGM).

**Theorem 4 (Lagrangian-IGM).**
The Lagrange expression for WLS is given by:

$$y = \left(\sum_{i=1}^{n}\sum_{j=1}^{n}(w_i - a_{ij}w_j)^2\right) + \dot{\lambda}\left(\left(\sum_{i=1}^{n} w_i\right) - 1\right). \tag{59}$$

The partial derivatives of $y$ with respect to $w_k$ and $\lambda$ are given by:

$$\frac{\partial y}{\partial w_k} = \left((n-1) + \sum_{i\neq k} a_{ik}^2\right) w_k + \left(\sum_{i\neq k}(-a_{ik} - a_{ki})w_i\right) + \lambda = 0, \qquad k = 1, \dots, n\,; \tag{60}$$

$$\frac{\partial y}{\partial \lambda} = \sum_{i=1}^{n} w_i - 1 = 0 \;\Rightarrow\; \sum_{i=1}^{n} w_i = 1. \tag{61}$$

where $\dot{\lambda} = 2\lambda$ is a Lagrange multiplier. The matrix expression of the system of linear equations above is

$$\ddot{G}\ddot{w}^T = \ddot{q}^T. \tag{62}$$

$\ddot{q}^T = [\mathbf{0}^T, 1]^T = [0,0,\dots,0,1]^T \in \mathbb{R}^{(n+1)\times 1}$, where $0 = [0,\dots,0] \in \mathbb{R}^n$ is a zero vector of length $n$. $\ddot{w}^T = [w, \lambda]^T = [w_1, \dots, w_n, \lambda]^T$ is a vector of solution weights and Lagrange multiplier. $\ddot{G}$ is also the symmetric augmented gram matrix. The 2x2 partition matrix $\ddot{G}$ is partitioned based on $\bar{G}$, i.e.,

$$\ddot{G} = \begin{vmatrix} \bar{G} & \boldsymbol{e^T} \\ e & 0 \end{vmatrix}. \tag{63}$$

The solution of Lagrangian Inverse Gram Matrix is given by:

$$\ddot{w} = (w, \lambda) = \ddot{q}\ddot{G}^{-1}. \tag{64}$$

**Proof.**

Unlike the constraint function setting in Eq.(51) presented by Chu et al.[19], the Lagrange expression based on objective and constraint functions in Eqs. (6) and (7) is of the form shown in Eq.(59). $\lambda$ is the Lagrange multiplier for the constraint function, $\sum_{i=1}^{n} w_i - 1 = 0$. To arrange the form in Eq.(59),

$$y = \left(\sum_{i=1}^{n}\sum_{j=1}^{n}\left(w_i^2 - 2a_{ij}w_jw_i + \left(a_{ij}w_j\right)^2\right)\right) - \dot{\lambda} + \dot{\lambda}\sum_{i=1}^{n} w_i. \tag{65}$$

When $i = j$, $a_{ij} = 1$, and thus

$$w_i^2 - 2a_{ij}w_jw_i + \left(a_{ij}w_j\right)^2 = 0. \tag{66}$$

Eq. (65) becomes:

$$y = \left(\sum_{i=1}^{n}\sum_{j\neq i}\left(w_i^2 - 2a_{ij}w_jw_i + \left(a_{ij}w_j\right)^2\right)\right) - \dot{\lambda} + \dot{\lambda}\sum_{i=1}^{n} w_i. \tag{67}$$

Taking the partial derivative of $y$ with respect to $w_k$, $k = 1, \dots, n$:

$$\frac{\partial y}{\partial w_k} = 2\sum_{j\neq k} w_k - 2\sum_{j\neq k} a_{kj}w_j - 2\sum_{i\neq k} a_{ik}w_i + 2\sum_{i\neq k}\left(a_{ik}^2 w_k\right) + \dot{\lambda} = 0. \tag{68}$$

To simplify the coefficient 2 in the form above, substitute

$$\lambda = \frac{\dot{\lambda}}{2} \tag{69}$$

into Eq. (68) and arrange terms:

$$\frac{\partial y}{\partial w_k} = \left((n-1)w_k + w_k \sum_{i \neq k} a_{ik}^2\right) + \left(-\sum_{i \neq k,i} a_{ik} w_i - \sum_{j \neq k,j} a_{kj} w_j\right) + \lambda = 0. \quad (70)$$

Finally, Eq. (60) holds. Taking partial derivative of $y$ with respect to $\lambda$, Eq. (61) holds. The Gram matrix $\ddot{G}$ of size $(n+1) \times (n+1)$ with coefficient elements is shown follows:

$$\ddot{G} = [\ddot{g}_{ij}] = \begin{bmatrix} \bar{G} & e^T \\ e & 0 \end{bmatrix} = \begin{bmatrix} (n-1) + \sum_{i \neq 1} a_{i1}^2 & -a_{21} - a_{12} & \cdots & -a_{n1} - a_{1n} & 1 \\ -a_{21} - a_{12} & (n-1) + \sum_{i \neq 2} a_{i2}^2 & \cdots & -a_{2n} - a_{n2} & 1 \\ \vdots & \vdots & \ddots & \vdots & \vdots \\ -a_{n1} - a_{1n} & -a_{2n} - a_{n2} & \cdots & (n-1) + \sum_{i \neq 3} a_{in}^2 & 1 \\ 1 & 1 & 1 & 1 & 0 \end{bmatrix} \quad (71)$$

The block $\bar{G}$ is proved in Eq. (46) Theorem 3. $\bar{G}$ is a block in $\ddot{G}$, which is non-singular and invertible. Thus, Eq. (63) holds. Eq, (62) holds and leads to the inverse solution shown in Eq.(64) . QED.

$\ddot{G}$ is invertible and satisfies the structural properties of Gram Matrix. The major difference between $\ddot{G}$ and $\bar{G}$ is the inclusion of the Lagrange multiplier $\lambda$ which makes $\ddot{G}$ inhomogeneous and non-singular when $A$ is perfectly consistent. Although several ways such as Gaussian elimination can solve the linear system, for the focus of this study, the inverse method is chosen as the solution shown in Eq.(64). $\tilde{G}_r = \bar{G} + rE$ may be used instead of $\bar{G}$ in $\tilde{G}$, i.e.,

$$\ddot{G}_r = \begin{bmatrix} \bar{G} + rE & e^T \\ e & 0 \end{bmatrix}. \quad (72)$$

Although Lagrange multiplier element in the last row and the other elements of $\ddot{G}^{-1}$ are changed, $w$ located in the last row of $\ddot{G}^{-1}$, which is the primary focus of interest, is the same, no matter what real scalar $r$ is used for $\tilde{G}$. To simply the calculation, $r = 0$ is used for Eq.(72) leading to Eq. (71), although there is nothing wrong with choosing another real value including $r = 1$ for $r$. In addition, the relations for $\bar{G}$, $\tilde{G}$, and $G$ used in $\ddot{G} = \ddot{G}_{r=0}$ is shown as follows.

$$\ddot{G} = \begin{vmatrix} \bar{G} & \boldsymbol{e^T} \\ e & 0 \end{vmatrix} = \begin{bmatrix} \tilde{G} - rE & e^T \\ e & 0 \end{bmatrix} = \begin{vmatrix} G - E & \boldsymbol{e^T} \\ e & 0 \end{vmatrix}. \quad (73)$$

### 3.4. Global Optimality and Uniqueness

To establish the theoretical validity of the Inverse Gram Matrix (IGM) framework, it is necessary to prove that the closed-form weight vector derived by IGM methods is not merely a stationary point, but the unique global minimizer of the Weighted Least Squares (WLS) optimization problem defined in Eqs. (6) and (7). In this section, Lemma 2 derives the exact Hessian matrix of the WLS objective function, and Theorem 5 provides the formal proof of strict convexity and unique global optimality.

**Lemma 2 (Hessian Matrix of the WLS Objective Function)**

Let $f(w) = \sum_{i=1}^{n}\sum_{j=1}^{n}\left(w_i - a_{ij}w_j\right)^2$ be the Weighted Least Squares (WLS) objective function defined in Eq. (6), and let $\overline{D} \in \mathbb{R}^{n(n-1)\times n}$ be the objective design matrix constructed by Algorithm 1 in Eq. (14). The objective function can be written in quadratic form as $f(w) = \left\|\overline{D}w^T\right\|^2 = w\left(\overline{D}^T\overline{D}\right)w^T$, and its Hessian matrix with respect to $w = [w_1, w_2, \dots, w_n]$ is:

$$\nabla^2 f(w) = 2\overline{D}^T\overline{D} \tag{74}$$

**Proof.**

By Algorithm 1 (Eq. (14)), each row $k \in \{1, \dots, n(n-1)\}$ of $\overline{D}$ represents the coefficient vector of a pairwise weighted differential relation $w_i - a_{ij}w_j = 0$ for distinct $i \neq j$. Stacking all $n(n-1)$ rows yields the vector $\overline{D}w^T$, so that:

$$f(w) = \sum_{i=1}^{n}\sum_{j=1}^{n}\left(w_i - a_{ij}w_j\right)^2 = \left\|\overline{D}w^T\right\|^2 = \left(\overline{D}w^T\right)^T\left(\overline{D}w^T\right) = w\left(\overline{D}^T\overline{D}\right)w^T \tag{75}$$

Let $M = \overline{D}^T\overline{D} \in \mathbb{R}^{n\times n}$. Since $M^T = \left(\overline{D}^T\overline{D}\right)^T = \overline{D}^T\overline{D} = M$, $M$ is a real symmetric matrix. Expanding the quadratic form component-wise gives

$$f(w) = \sum_{i=1}^{n}\sum_{j=1}^{n} M_{ij}\, w_i w_j \tag{76}$$

Taking the first-order partial derivative with respect to $w_k$ ($k = 1, \dots, n$):

$$\frac{\partial f(w)}{\partial w_k} = \frac{\partial}{\partial w_k}\left(\sum_{i=1}^{n}\sum_{j=1}^{n} M_{ij}\, w_i w_j\right) = 2\sum_{j=1}^{n} M_{kj}\, w_j \tag{77}$$

Taking the second-order partial derivative with respect to $w_l$ ($l = 1, \dots, n$):

$$\frac{\partial^2 f(w)}{\partial w_k \partial w_l} = \frac{\partial}{\partial w_l}\left(2\sum_{j=1}^{n} M_{kj}\, w_j\right) = 2M_{kl} = 2\left(\overline{D}^T\overline{D}\right)_{kl} \tag{78}$$

For all entry pairs $(k, l) \in \{1, \dots, n\}^2$, Eq.(74) is obtained from Eq. (78). QED.

Having established in Lemma 2 that the Hessian matrix of the Weighted Least Squares (WLS) objective function is $\nabla^2 f(w) = 2\overline{D}^T\overline{D}$, the quadratic form of this Hessian is evaluated on the tangent space of the feasible weight simplex $\mathcal{S} = \{w \in \mathbb{R}^n : ew^T = 1, w > \mathbf{0}\}$ in Theorem 5. Theorem 5 demonstrates strict convexity and integrates this result with the Lagrangian stationary system from Theorem 4 to prove that the Inverse Gram Matrix (IGM) solution is the unique global minimizer.

**Theorem 5 (Global Optimality and Uniqueness of IGM for WLS Solution).**

The priority weight vector $w^* = [w_1^*, w_2^*, \dots, w_n^*]$ derived by the Inverse Gram Matrix (IGM) methods is the unique global minimizer of the Weighted Least Squares (WLS) optimization problem defined in Eqs. (6) and (7).

**Proof.**

Let $D = \begin{bmatrix} \overline{D} \\ e \end{bmatrix}$ be the design matrix constructed by Algorithm 1 in Eq. (14), where $\overline{D} \in \mathbb{R}^{n(n-1)\times n}$ contains the pairwise differential coefficients $w_i - a_{ij} w_j = 0$, and $e = [1,1,\dots,1] \in \mathbb{R}^{1\times n}$.

**1.** Strict Convexity of the Objective Function (via Lemma 2 and Lemma 1):

By Lemma 2, the WLS objective function is $f(w) = \|\overline{D} w^T\|^2 = w\left(\overline{D}^T \overline{D}\right) w^T$, and its Hessian matrix is $\nabla^2 f(w) = 2\overline{D}^T \overline{D}$.

For any non-zero displacement vector $h = w^{(1)} - w^{(2)}$ between two distinct feasible vectors

$$w^{(1)}, w^{(2)} \in \mathcal{S} = \{w \in \mathbb{R}^n : ew^T = 1, w > \mathbf{0}\}(h \neq \mathbf{0}), \tag{79}$$

the normalization constraint requires:

$$eh^T = e\left(w^{(1)}\right)^T - e\left(w^{(2)}\right)^T = 1 - 1 = 0 \Longrightarrow h \in \mathrm{null}(e) \tag{80}$$

By Lemma 1, $\mathrm{null}(D) = \mathrm{null}(\overline{D}) \cap \mathrm{null}(e) = \{\mathbf{0}\}$. Since $h \in \mathrm{null}(e)$ and $h \neq \mathbf{0}$, $h$ cannot lie in $\mathrm{null}(\overline{D})$, which implies $\overline{D} h^T \neq \mathbf{0}$. Therefore:

$$h\left(\overline{D}^T \overline{D}\right) h^T = \|\overline{D} h^T\|^2 > 0 \tag{81}$$

This proves that the Hessian $\nabla^2 f(w) = 2\overline{D}^T \overline{D}$ is strictly positive-definite on the tangent space of the constraint hyperplane $\{h: eh^T = 0\}$. Consequently, the objective function $f(w)$ is strictly convex on the feasible set $\mathcal{S}$.

2. Integration with Theorem 4 (Uniqueness & Global Optimality):

By Theorem 4, setting the gradient of the Lagrangian $y(w, \lambda)$ in Eq. (59) to zero yields the first-order necessary conditions (KKT system) represented by the block linear system $\overline{G}[w, \lambda]^T = q^T$ (Eq. (62)). For a strictly convex objective function defined on a convex feasible set, first-order stationary conditions are both necessary and sufficient for global optimality.

Furthermore, by Lemma 1, the Gram matrix $G = D^T D = \overline{D}^T \overline{D} + e^T e$ is strictly invertible ($\det(G) > 0$) and $\overline{G}$ is non-singular, guaranteeing that the stationary point $[w^*, \lambda^*] = q\overline{G}^{-1}$ derived in Theorem 4 (Eq. (64)) is unique and achieves the strict global minimum of the WLS objective function.

QED.

**Remark 2 (Equivalence of IGM Variations):**

Although Pseudo-IGM (PIGM), Normalized-IGM (NIGM), and Lagrangian-IGM (LIGM) are constructed via different matrix operations, they are algebraically equivalent formulations of the unique global WLS minimizer $w^*$. LIGM provides the explicit KKT multiplier derivation, PIGM provides the design matrix pseudoinverse formulation, and N-IGM provides the element-form Gram matrix representation.

Table 2. Summary of a family of three Inverse Gram Matrix Methods

| Name | Pseudo-IGM (PIGM) | Normalized-IGM (NIGM) | Lagrangian -IGM (LIGM) |
|---|---|---|---|
| Form | $w = \frac{eG^{-1}}{eG^{-1}e^T} = \frac{e(D^TD)^{-1}}{e(D^TD)^{-1}e^T}$ | $w = \frac{e\tilde{G}^{-1}}{e\tilde{G}^{-1}e^T} = \frac{e(\bar{G}+r)^{-1}}{e(\bar{G}+r)^{-1}e^T}$ | $\ddot{w} = (w, \lambda) = \ddot{q}\ddot{G}^{-1}$ |
| Gram matrix | $g_{ij} = \begin{cases} (n-1) + \sum_k a_{kj}^2 & i = j \\ 1 - a_{ij} - a_{ji} & i \neq j \end{cases}$ | $\tilde{g}_{ij} = r + \bar{g}_{ij}$ $= \begin{cases} r + (n-1) + \sum_{k \neq i} a_{kj}^2 & i = j \\ r - a_{ij} - a_{ji} & i \neq j \end{cases}$ | $\ddot{G} = \begin{vmatrix} \bar{G} & \boldsymbol{e^T} \\ e & 0 \end{vmatrix}$ $\bar{G} = C = B$ |
| Equation reference number | Eqs. (23) , (14) and (33) (or (34)). | Eqs. (42)-**Error! Reference source not found.** or Eq. (49) | Eqs. (63) and (64) |
| Matrix Dimension | $n \times n$ | $n \times n$ | $(n+1) \times (n+1)$ |
| Complexity | $\mathcal{O}(n^3)$ | $\mathcal{O}(n^3)$ | $\mathcal{O}(n^3)$ |
| Practitioner Selection Criteria | Best for standard pseudoinverse system formulations where design matrix $D$ is explicitly formed. | Recommended for direct code implementation. Bypasses design matrix $D$; parameter $r = 1$ guarantees invertibility for consistent PRMs. | Recommended for exact KKT analysis. Solves $w$and Lagrange multiplier $\lambda$ simultaneously without requiring post-hoc normalization. |

**Remark 3 (Computational Complexity and Practitioner Selection Criteria):**

While all three proposed IGM variants—Pseudo-IGM (PIGM), Normalized-IGM (N-IGM), and Lagrangian-IGM (LIGM)—share the exact same asymptotic computational complexity class of $\mathcal{O}(n^3)$ operations for direct matrix inversion and yield mathematically identical priority weight vectors $w^*$, Table 2 summarizes their structural expressions, matrix dimensions, and practical selection guidelines:

1. Pseudo-IGM (PIGM): Evaluates $w = \frac{eG^{-1}}{eG^{-1}e^T}$ using the $n \times n$ Gram matrix $G = D^TD$. It is best suited for standard linear system and pseudoinverse matrix frameworks where the design matrix $D \in \mathbb{R}^{(n(n-1)+1)\times n}$ is explicitly generated.

2. Normalized-IGM (N-IGM): Evaluates $w = \frac{e(G+rE)^{-1}}{e(G+rE)^{-1}e^T}$ using closed-form element formulas $g_{ij}$. It is recommended for direct software implementation, as it completely bypasses the construction of the large design matrix $D$ and guarantees non-singular matrix invertibility via regularization parameter $r \neq 0$ (e.g., $r = 1$) even for perfectly consistent PRMs.

3. Lagrangian-IGM (LIGM): Solves the partitioned $(n+1) \times (n+1)$ block KKT linear system $[w, \lambda] = qG^{-1}$. It is recommended for analytical optimization tracking, as it computes both the priority weight vector $w$ and the exact Lagrange multiplier $\lambda$ in a single matrix inversion step without requiring post-hoc vector normalization.

**Remark 4 (Ratio Scale Preservation under Inconsistency).**

In AHP, priority vectors $w = [w_1, \ldots, w_n]^T$ represent a ratio scale on the unit simplex $\sum_{i=1}^{n} w_i = 1$. When a PRM is inconsistent ($a_{ij} \neq w_i/w_j$), IGM methods preserve the ratio scale property by deriving the exact analytical global minimizer of the WLS objective $\min \sum_{i,j} \left(w_i - a_{ij} w_j\right)^2$. By minimizing weighted ratio deviations in closed form without numerical solver approximation error, IGM preserves maximum ratio scale fidelity to the expert's judgment distribution.

## 4. Numerical Examples

One consistent and one inconsistent PRMs are used to demonstrate the calculation steps of three proposed Inverse Gram Matrix Methods summarized in Table 2 in Examples 1 and 2, respectively. For the baseline testing, when a PRM is consistent, $w$ has a unique solution, which must be produced by a valid prioritization method. To solve the optimization model in Eqs. (6) and (7), the NLOPT_GN_ISRES solver algorithm of nloptr package [40] in R is used, whilst the *solve()* function in R is used to find the inverse of a matrix. For the following calculations, values are rounded to three decimal places for display, but full internal precision is retained during computation. When priority weights $w$ are rounded to three decimal places, the WLS optimization model of Eqs.(6) and (7) using optimization software produce the same results as IGM methods below.

To facilitate practical implementation and computational reproducibility, all numerical calculations demonstrated in this section are fully implemented in the open-source R package *igm4ahp* [36]. Complete reproducible R code and step-by-step computational workflows for the examples are available in the package vignette documentation at https://kkfyuen.github.io/igm4ahp/articles/section4_igm_examples.html.

### 4.1. Example 1

For the baseline of validity, a prioritization method must produce the unique priority solution for a perfectly consistent PRM. Given a perfectly consistent $4 \times 4$ PRM below:

$$A_1 = \begin{bmatrix} 1 & 2 & 4 & 8 \\ \frac{1}{2} & 1 & 2 & 4 \\ \frac{1}{4} & \frac{1}{2} & 1 & 2 \\ \frac{1}{8} & \frac{1}{4} & \frac{1}{2} & 1 \end{bmatrix}.$$

By Algorithm 1, a $(12+1)\times 4$ design matrix $D$ is of the form below:

$$D=\begin{bmatrix}\bar{D}\\ e\end{bmatrix}=\begin{bmatrix} & \bar{D} & & \\ 1 & 1 & 1 & 1\end{bmatrix}, \bar{D}=\begin{bmatrix}1 & -2 & 0 & 0\\ -\frac{1}{2} & 1 & 0 & 0\\ 1 & 0 & -4 & 0\\ -\frac{1}{4} & 0 & 1 & 0\\ 1 & 0 & 0 & -8\\ -\frac{1}{8} & 0 & 0 & 1\\ 0 & 1 & -2 & 0\\ 0 & -\frac{1}{2} & 1 & 0\\ 0 & 1 & 0 & -4\\ 0 & -\frac{1}{4} & 0 & 1\\ 0 & 0 & 1 & -2\\ 0 & 0 & -\frac{1}{2} & 1\end{bmatrix}.$$

The last row of $D$ is $e$ with respect to the coefficients of the weight sum equal to unity function and theremaining rows of $D$ form the $12\times 4$ objective differential design matrix, $D$. The $4\times 4$ gram matrix computed by matrix multiplication shown in Eq.(15) is of the form below.

$$G=D^TD=\begin{bmatrix}4.328 & -1.5 & -3.25 & -7.125\\ -1.5 & 8.312 & -1.5 & -3.250\\ -3.25 & -1.5 & 24.25 & -1.5\\ -7.125 & -3.25 & -1.5 & 88\end{bmatrix}.$$

Alternatively, without constructing a $D$ matrix, $G$ can be computed in element form via Eq. (33). For example, two entries of $G$ are computed as below.

$$g_{11}=(4-1)+\left(1+\left(\frac{1}{2}\right)^2+\left(\frac{1}{4}\right)^2+\left(\frac{1}{8}\right)^2\right)=4.328;$$

$$g_{12}=1-2-\frac{1}{2}=-1.5.$$

The inverse gram matrix computed using the *solve()* function in R is of the form below.

$$G^{-1}=\begin{bmatrix}0.357 & 0.087 & 0.055 & 0.033\\ 0.087 & 0.145 & 0.021 & 0.013\\ 0.055 & 0.021 & 0.050 & 0.006\\ 0.033 & 0.013 & 0.006 & 0.015\end{bmatrix}.$$

By Pseudo-IGM, summing each row or column of the symmetric $G^{-1}$ to have

$$v=eG^{-1}=(0.533,0.267,0.133,0.067);$$

$$eG^{-1}e^T=\sum v=1;$$

$$w = \frac{eG^{-1}}{eG^{-1}e^T} = (0.533,0.267,0.133,0.067).$$

For the method mentioned in *Blankmeyer* [39], by Eq.(41),

$$\bar{G} = C = \bar{D}^T\bar{D} = G - 1 = \begin{bmatrix} 3.328 & -2.5 & -4.25 & -8.125 \\ -2.5 & 7.312 & -2.5 & -4.250 \\ -4.25 & -2.5 & 23.25 & -2.5 \\ -8.125 & -4.25 & -2.5 & 87 \end{bmatrix}.$$

For a perfectly consistent *A*, as $\bar{G}$ (or $C$) is a singular matrix ($\det \bar{G} = 0$) and not invertible, Eqs. (8)-(10) cannot produce the *w* solution. NIGM proposed in Theorem 3 is used to solve this problem.

$$w = \frac{e\tilde{G}^{-1}}{e\tilde{G}^{-1}e^T} = \frac{e(\bar{G}+r)^{-1}}{e(\bar{G}+r)^{-1}e^T}.$$

*r* is a non-zero real number ($r \neq 0$).

For *r* = 1, the result is the same as PIGM. For *r* = 5 as an illustration of NIGM, i.e.,

$$\tilde{G} = \begin{bmatrix} \bar{D} \\ e\sqrt{5} \end{bmatrix}^T \begin{bmatrix} \bar{D} \\ e\sqrt{5} \end{bmatrix} = \bar{G} + 5E .$$

The inverse of $\tilde{G}$ is:

$$\tilde{G}^{-1} = \begin{bmatrix} 0.1300 & -0.0263 & 0.0015 & 0.0047 \\ -0.0263 & 0.0881 & 0.0070 & 0.0014 \\ 0.0015 & 0.0070 & 0.0361 & -0.0010 \\ 0.0047 & 0.0014 & -0.0010 & 0.0111 \end{bmatrix}.$$

Finally,

$$w = \frac{e\tilde{G}^{-1}}{e\tilde{G}^{-1}e^T} = \frac{(0.1067,0.0533,0.0267,0.0133)}{0.2} = (0.533,0.267,0.133,0.067).$$

Similarly, if we take any other non-zero value for *r*, *w* is the same as PIGM after the normalization procedure above.

By Lagrangian-IGM in Theorem 4, the augmented Lagrange gram matrix is

$$\ddot{G} = \begin{vmatrix} \bar{G} & \boldsymbol{e^T} \\ e & 0 \end{vmatrix} = \begin{bmatrix} 3.328 & -2.5 & -4.25 & -8.125 & 1 \\ -2.5 & 7.312 & -2.5 & -4.250 & 1 \\ -4.25 & -2.5 & 23.25 & -2.5 & 1 \\ -8.125 & -4.25 & -2.5 & 87 & 1 \\ 1 & 1 & 1 & 1 & 0 \end{bmatrix}.$$

Its inverse matrix is:

$$\ddot{G}^{-1} = \begin{bmatrix} 0.0730 & -0.0547 & -0.0158 & -0.0024 & 0.5333 \\ -0.0548 & 0.0738 & -0.0141 & -0.0050 & 0.2667 \\ -0.0158 & -0.0141 & 0.0326 & -0.0028 & 0.1333 \\ -0.0024 & -0.0050 & -0.0028 & 0.0102 & 0.0666 \\ 0.5333 & 0.2667 & 0.1333 & 0.0666 & 0 \end{bmatrix}.$$

Finally, the last row or column of $\ddot{G}^{-1}$ excluding the last element of Lagrange Multiplier is the solution of $w$, i.e.,

$$\ddot{w} = [w, \lambda] = \ddot{q}\ddot{G}^{-1} = (0.533,0.267,0.133,0.067,0).$$

## 4.2. Example 2

The inconsistent PRM below evaluates six criteria listed in sequence to select a high school: learning, friends, school life, vocational training, college preparation, and music classes [5].

$$A_2 = \begin{bmatrix} 1 & 4 & 3 & 1 & 3 & 4 \\ 1/4 & 1 & 7 & 3 & 1/5 & 1 \\ 1/3 & 1/7 & 1 & 1/5 & 1/5 & 1/6 \\ 1 & 1/3 & 5 & 1 & 1 & 1/3 \\ 1/3 & 5 & 5 & 1 & 1 & 3 \\ 1/4 & 1 & 6 & 3 & 1/3 & 1 \end{bmatrix}.$$

Consistency index and ratio are 0.3 and 0.24 respectively. $G$ can be computed by element form with Eq. (33) or Eq. (34) and Eq.(15).

$$G = \begin{bmatrix} 7.347 & -3.25 & -2.333 & -1 & -2.333 & -3.25 \\ -3.25 & 48.132 & -6.143 & -2.333 & -4.2 & -1 \\ -2.333 & -6.143 & 150 & -4.2 & -4.2 & -5.167 \\ -1 & -2.333 & -4.2 & 26.04 & -1 & -2.333 \\ -2.333 & -4.2 & -4.2 & -1 & 16.191 & -2.333 \\ -3.25 & -1 & -5.167 & -2.333 & -2.333 & 32.139 \end{bmatrix}.$$

The inverse gram matrix is of the form below.

$$G^{-1} = \begin{bmatrix} 0.166 & 0.016 & 0.005 & 0.012 & 0.033 & 0.021 \\ 0.016 & 0.023 & 0.002 & 0.004 & 0.009 & 0.004 \\ 0.005 & 0.002 & 0.007 & 0.002 & 0.003 & 0.002 \\ 0.012 & 0.004 & 0.002 & 0.04 & 0.006 & 0.005 \\ 0.033 & 0.009 & 0.003 & 0.006 & 0.072 & 0.01 \\ 0.021 & 0.004 & 0.002 & 0.005 & 0.01 & 0.035 \end{bmatrix}.$$

By PIGM, the weights are

$$w = \frac{eG^{-1}}{eG^{-1}e^T} = \frac{(0.254, 0.057, 0.021, 0.069, 0.134, 0.077)}{(0.254 + 0.057 + 0.021 + 0.069 + 0.134 + 0.077)}$$
$$= (0.415, 0.094, 0.035, 0.112, 0.219, 0.125)\ .$$

By NIGM, the easiest way is to choose *r*=0 since $A_2$ is not perfectly consistent. As

$$\tilde{G} = \begin{bmatrix} \overline{D} \\ 0 \end{bmatrix}^T \begin{bmatrix} \overline{D} \\ 0 \end{bmatrix} = \bar{G} = G - 1,$$

the inverse matrix is

$$\tilde{G}^{-1} = \bar{G}^{-1} = \begin{bmatrix} 0.333 & 0.053 & 0.019 & 0.057 & 0.121 & 0.072 \\ 0.053 & 0.032 & 0.005 & 0.014 & 0.029 & 0.015 \\ 0.019 & 0.005 & 0.008 & 0.006 & 0.011 & 0.006 \\ 0.057 & 0.014 & 0.006 & 0.052 & 0.03 & 0.019 \\ 0.121 & 0.029 & 0.011 & 0.03 & 0.118 & 0.036 \\ 0.072 & 0.015 & 0.006 & 0.019 & 0.036 & 0.05 \end{bmatrix}.$$

Taking the sum of rows or columns of the matrix above:

$$e\tilde{G}^{-1} = (0.655\ 0.148\ 0.055\ 0.177\ 0.346\ 0.198).$$

By Eq.(42), $w$ is the same as PIGM. Please note that if PRM is perfectly consistent, $r=0$ leads to error due to non-invertible gram matrix. This is shown in the previous example. For the easiest and safest calculation, $r$=1 is recommended, under which NIGM is identical to PIGM.

By Lagrangian-IGM, the Lagrange gram matrix based on Eq. (72) where $\tilde{r} = 1$ and $G = \bar{G} + 1$ is shown as below.

$$\ddot{G} = \begin{vmatrix} G & \boldsymbol{e^T} \\ e & 0 \end{vmatrix}.$$

The Lagrange inverse gram matrix is

$$\ddot{G}^{-1} = \begin{bmatrix} 0.061 & -0.008 & -0.004 & -0.017 & -0.022 & -0.01 & 0.415 \\ -0.008 & 0.018 & 0 & -0.003 & -0.003 & -0.004 & 0.094 \\ -0.004 & 0 & 0.006 & -0.001 & -0.001 & -0.001 & 0.035 \\ -0.017 & -0.003 & -0.001 & 0.032 & -0.009 & -0.004 & 0.112 \\ -0.022 & -0.003 & -0.001 & -0.009 & 0.042 & -0.007 & 0.219 \\ -0.01 & -0.004 & -0.001 & -0.004 & -0.007 & 0.025 & 0.125 \\ 0.415 & 0.094 & 0.035 & 0.112 & 0.219 & 0.125 & -1.633 \end{bmatrix}.$$

By Eq. (64), the solution for $w, \lambda$ is the last row or column of $\ddot{G}^{-1}$. No matter what real value is chosen for $r$, $w$ is the same as $\lambda$ is changed accordingly. A small value like $r = 1$ or $r = 0$ is recommended to avoid unnecessary extra calculation.

## 5. Simulations and Verifications

As mentioned earlier, many comparisons and discussions regarding prioritization methods, including WLS, exist in the literature, which is not the primary focus of this paper. The objective of this paper is to propose IGM methods, which are interestingly the exact closed-form solution to WLS. Therefore, IGM methods produce the least WLS error over all other prioritization methods. This paper highlights the computational efficiency and accuracy of the inverse matrix methods in deriving exact values, instead of relying on numerical optimization solvers for WLS to find approximate values. To further verify the validity of solutions from the three proposed IGM methods summarized in Table 2—and specifically to investigate whether any edge or special cases were overlooked during algorithm design and implementation—the pseudocode of the proposed verification algorithm is presented in Algorithm 2. The 9-point scale is used for $\tau$. The maximum number of criteria, $n_{\max}$, is set to 15. Benchmarking was performed on an HP OMEN 16 Laptop equipped with an Intel i9-13900HX CPU and 32 GB RAM, using a single CPU core process for consistent timing comparisons.

Let $w^{WLS}$, $w^{PIGM}$, $w^{NIGM(r)}$, $w^{LIGM}$ and $w^{LIGM(r)}$ be priority weight vectors obtained by WLS (Eqs. (6) and (7)), PIGM (Eq. (23)), NIGM($r$) (Eq. (42)), LIGM (Eq.(64)) and LIGM($r$) (Eq.(72)), respectively. Although $r$ is a non-zero real number, it is set between -1000 and 1000 for verification, as a large value for $r$ is not recommended due to computational workload. In principle, the following relationships are verified.

$$errorDectected = \begin{cases} False & w^{WLS} = w^{PIGM} = w^{NIGM(r)} = w^{LIGM} = w^{LIGM(r)} \\ True & Otherwise \end{cases} \tag{82}$$

In practise, the error detection is implemented as two parts due to different computational efficiency and floating-point comparisons. $testwlsOpt$ in Algorithm 2 is used to set which part is tested.

Simulation 1 tests the numerical equality $w^{PIGM} = w^{NIGM(r)} = w^{LIGM} = w^{LIGM(r)}$. The main procedures are shown in Step 2b of Algorithm 2. PIGM, LIGM, and NIGM theoretically produce the same results, which however are not always identical in raw computation. Executing floating-point operations may induce rounding and truncation errors, where a very small difference could trigger false equality failures. Therefore, the IGM equality verification in Eq.(82) is implemented via Eq. (84). The inverse matrix function is performed by the solve() function in R. The simulation of one million random PRMs calculated by the four IGM methods took only 21.94 minutes (1,316.48 seconds). The simulation results are available in Supplementary Data 1. Complete reproducible R code and step-by-step computational workflows for Simulation 1 are available in the package vignette documentation at https://kkfyuen.github.io/igm4ahp/articles/section5_simulation1_igm.html. No error was found across the one million random instances, confirming that all four methods produce identical results. The rounding decimal places parameter, $\varepsilon$, is set to 8 for the round() function.

Simulation 2 evaluates the equality between the numerical WLS optimization solution and the proposed Inverse Gram Matrix (IGM) methods (specifically Pseudo-IGM in this setup), as detailed in Step 2c of Algorithm 2. Because numerical optimization solvers generally obtain approximate solutions that converge toward the exact solution, a lower approximation error for $w^{WLS}$ can typically be achieved by setting a higher computational search budget. Therefore, the WLS and IGM equality verification in Eq.(82) is implemented via Eq.(85) after Eq.(84) has been verified. The rounding precision parameter, $\varepsilon$, is set to 4 decimal places for the rounding function.

The NLOPT_GN_ISRES solver from the nloptr package [40] in R is employed to solve the WLS model in Eqs. (6) and (7), initialized with the PIGM weight vector. To ensure reliable global convergence without premature termination, the maximum number of function evaluations is set to 500,000 for the fractional stopping tolerance of $10^{-6}$. Because NLOPT_GN_ISRES is a stochastic evolutionary algorithm operating on matrix dimensions up to $n = 15$, smaller evaluation budgets (e.g., $< 100{,}000$) risk terminating in intermediate generations prior to full convergence on higher-dimensional simplexes. Benchmarking confirms that 500,000 evaluations provide sufficient generational depth for $n \leq 15$ to reach $|\, f(w_{\mathrm{WLS}}) - f(w_{\mathrm{IGM}})\, | < 10^{-6}$, ensuring that no false errors occur due to incomplete numerical optimization. If an error were detected (i.e., weights deviating with higher WLS objective error), the computational budget would be increased and the simulation re-run.

While iterative non-linear optimization solvers scale with heavy evaluation workloads ($\mathcal{O}(K \cdot n^2)$, where each of the $K$ function evaluations requires $\mathcal{O}(n^2)$ operations across $n^2$ pairwise terms), closed-form IGM methods rely on direct matrix inversion requiring only $\mathcal{O}(n^3)$ operations. Comparing the two complexity classes yields an Efficiency Ratio:

$$\frac{\text{Iterative Solver Cost}}{\text{IGM Closed-Form Cost}} = \frac{\mathcal{O}(K \cdot n^2)}{\mathcal{O}(n^3)} = \mathcal{O}\left(\frac{K}{n}\right) \tag{83}$$

Because global evolutionary solvers require a high evaluation budget ($K = 500{,}000 \gg n$) to ensure global convergence without premature termination across high-dimensional simplexes ($n \leq 15$), the efficiency ratio $K/n \approx 33{,}333$ reveals that the numerical solver requires over 33,000 times more computational work per matrix instance. In practice, the iterative optimization approach requires over $10^8$ operations per matrix (averaging $\sim 15.7$ seconds), whereas direct IGM inversion requires only $\sim 3.3 \times 10^3$ operations (averaging $< 0.0013$ seconds), establishing decisive computational efficiency for the proposed closed-form methods.

For simulation 2, 10,000 random PRMs were evaluated over 156,882.8 seconds (43.58 hours). The simulation results are available in Supplementary data 2. Complete reproducible R code and step-by-step computational workflows for Simulation 2 are available in the package vignette documentation at https://kkfyuen.github.io/igm4ahp/articles/section5_simulation2_wlsopt.html. No error or unexpected case was found across any of the generated random PRM instances. It can be concluded that IGM methods produce the exact closed-form solution for WLS of Eqs.(6) and (7), and the results can be used to verify the solutions of optimization solvers. By comparing the computational time in Simulations 1 and 2, obtaining higher precision from the numerical solver is very computationally expensive. The IGM methods can be used to quickly determine the weights with minimal computational effort.

**Remark 5 (Solver Termination Criteria, Precision Rationale, and Residual Distribution):**

In practical AHP decision models, priority weights are standardly reported to two decimal places. In Simulation 1, direct closed-form IGM evaluations were verified to eight decimal places ($\epsilon = 8$), yielding zero error across 1,000,000 random instances. In Simulation 2, the iterative numerical solver (NLOPT_GN_ISRES) was terminated using a maximum budget of 500,000 function evaluations and a fractional stopping tolerance of $10^{-6}$. A verification rounding threshold of $\epsilon = 4$ decimal places was selected to match this $10^{-6}$ solver tolerance. At $\epsilon = 4$, direct weight comparisons between numerical WLS and IGM yielded zero discrepancy across all 10,000 instances. Furthermore, across all 10,000 random PRMs in Simulation 2, the actual double-precision objective function residual $|f(w_{\text{WLS}}) - f(w_{\text{IGM}})|$ remained strictly bounded within $[0, 1.4 \times 10^{-7}]$, confirming that numerical optimization converges directly to the exact closed-form IGM solution. Because closed-form IGM matrix inversion eliminates numerical solver approximations, increasing solver search depth to higher decimal precision consumes unnecessary computational resources without changing the scientific conclusion.

Algorithm 2: **(**$Verification(N, n^{max}, \tau, testwlsOpt, \varepsilon)$)

**Input**: $N$: Testing sample size; $n^{max}$: maximum number of criteria; $\tau$: $z$-point scale for $a_{ij}$; $testwlsOpt$: False for Simulation 1; True to Simulation 2. $\varepsilon$: significant decimal places

**Step 1**: Initialize the values.

times = 1 # case number or counter

error = 0 # error checking

flag = TRUE # flag to continue or exist the while-loop.

**Step 2**: Perform comparisons for $N$ samples.

While (flag):

**Step 2a**: Generate a random PRM.

# Generate a random integer between 3 and $n^{max}$ for *n*.

$n = random(3, n^{max})$

# Generate a non-zero random number between -1000 and 1000.

$r = random(-1000, 1000, \neg 0)$

# Generate a random $n \times n$ PRM based $\tau$ scale.

$A = generatePRM(n, \tau)$

**Step 2b**: Determine to perform Simulation 1

**If (**$testwlsopt = FALSE$)

# Obtain *w* with four IGM methods.

$w^{PIGM} = PIGM(A)$

$w^{NIGM(r)} = NIGM(A, r)$

$w^{LIGM} = LIGM(A)$

$w^{LIGM(r)} = LIGM(A, r)$

# Test if any error.

$$error = round\left( \sum_{k \in \left\{ \substack{NIGM(r), \\ LIGM, LIGM(r)} \right\}} \left| round(w^{PIGM}, \varepsilon) - round\left(w^k, \varepsilon\right) \right|, (\varepsilon - 1) \right) \quad (84)$$

**Step 2c**: Determine to perform Simulation 2

**If (**$testwlsopt = TRUE$)

#Obtain *w* with difference methods.

$w^{WLS} = wlsOpt(A)$

$w^{PIGM} = PIGM(A)$

# Test if any error.

$$error = round(|round(w^{WLS}, \varepsilon) - round(w^{PIGM}, \varepsilon)|, \varepsilon - 1) \quad (85)$$

**Step 3c**: determine the exit flag for the while-loop.

# Increment the counter by 1

*times* += 1

# If *N* samples are performed or an error is found, exit the while-loop. Otherwise, continue.

$$flag = \begin{cases} \text{FALSE}, error \neq 0 \; or \; times > N \\ \text{TRUE}, \quad otherwise \end{cases}$$

End While

# If an error is found, print the relevant information for further study.

# Otherwise, no error is found, i.e., error = 0.

$$errorDectected = \begin{cases} True & error \neq 0 \\ FALSE & error = 0 \end{cases} \quad (86)$$

**Return**: $errorDectected$

## 6. Conclusions and Future Directions

The primary initiative of this study was to propose the Pseudo-Inverse Gram Matrix (PIGM) method to solve Pairwise Reciprocal Matrices (PRMs) in AHP. Interestingly, the PIGM method produces the exact closed-form solution for WLS without relying on numerical optimization solvers. The closed-form solutions of WLS from previous studies (e.g., [19, 39]) were incomplete, lacking formal proofs, step-by-step demonstrations, and simulation verifications. This study further proposes Normalized-IGM (NIGM) and Lagrangian-IGM (LIGM) to fill these research gaps. The closed-form solutions of WLS can thus be clearly explained from different mathematical perspectives. Lagrangian-IGM solves for the Lagrange multiplier directly without requiring post-hoc normalization, whilst Pseudo-IGM requires post-hoc normalization without needing to evaluate the Lagrange multiplier. Normalized-IGM provides a generalized formulation extending the properties of Pseudo-IGM, specifically explaining the element-form entries of the Gram matrix. Two numerical examples demonstrate the calculation steps, similarities, and distinctions among the three IGM methods. Extensive Monte Carlo simulations verify the numerical equivalence of the three IGM methods and confirm that they produce the exact global minimizer of WLS. To facilitate practical adoption and computational reproducibility, all proposed algorithms have been implemented in the open-source R package *igm4ahp* [36], providing researchers and practitioners with an accessible, high-performance toolkit for exact AHP prioritization. This paper provides computationally simple, closed-form alternatives for WLS prioritization in AHP, using direct matrix inversion instead of intensive numerical optimization.

Furthermore, while classic AHP literature discusses rank reversal phenomena when candidate sets are modified, relative measurement naturally reflects contextual shifts in judgment distributions upon introducing new alternatives. Because the proposed IGM methods furnish the exact closed-form solution for the Weighted Least Squares (WLS) objective on any static Pairwise Reciprocal Matrix, they preserve full fidelity to the WLS model without numerical solver error. Investigating rank preservation characteristics across varying prioritization norms (e.g., WLS, LLSM, and Eigenvector methods) under dynamic criteria and alternative additions remains a valuable direction for future research.

## Acknowledgements

Sincere thanks are extended to the editors and anonymous referees for their time and effort in improving this manuscript.

## Funding

This research was not been funded by any company or organization.

## Author Contributions

This research was solely conducted by the author, K.K.F. Yuen.

## Conflicts of Interest

The author declares no conflicts of interest.

## Ethical approval and consent to participate

Not applicable.

## Consent for publication

Not applicable.

## Data and Code Availability

- The package user manual, tutorials, and numerical demonstrations implementing the algorithms and methods in this paper are hosted at https://kkfyuen.github.io/igm4ahp/.
- The igm4ahp R package and simulation datasets (Supplementary Data 1 and 2) are available on GitHub at https://github.com/kkfyuen/igm4ahp.
- The complete code and data are permanently archived on Zenodo at https://doi.org/10.5281/zenodo.24223557 [36].
- The preprint of this work is available at https://doi.org/10.48550/arXiv.2401.01190.

## References

1. Faliszewski, P., E. Hemaspaandra, and L.A. Hemaspaandra, *Using complexity to protect elections.* Commun. ACM, 2010. **53**(11): p. 74–82.
2. Koczkodaj, W.W., et al., *Important Facts and Observations about Pairwise Comparisons.* Fundamenta Informaticae, 2016. **144**: p. 1-17.
3. Thurstone, L., *A law of comparative judgment.* Psychological Review, 1927. **34**(4): p. 273-286.
4. Saaty, T.L., *A scaling method for priorities in hierarchical structures.* Journal of Mathematical Psychology, 1977. **15**(3): p. 234-281.
5. Saaty, T.L., *Analytic Hierarchy Process: Planning, Priority, Setting, Resource Allocation.* McGraw-Hill, New York, 1980.
6. Saaty, T.L., *Theory and Applications of the Analytic Network Process: Decision Making with Benefits, Opportunities, Costs, and Risks.* RWS Publications, 2005.
7. Yuen, K.K.F., *Cognitive network process with fuzzy soft computing technique for collective decision aiding.* The Hong Kong Polytechnic University, 2009. **Ph.D. thesis**.
8. Yuen, K.K.F., *Pairwise opposite matrix and its cognitive prioritization operators: Comparisons with pairwise reciprocal matrix and analytic prioritization operators.* Journal of the Operational Research Society, 2012. **63**(3): p. 322-338.
9. Dyer, J.S., *Remarks on the analytic hierarchy process.* Management Science 1990. **36** (3): p. 249-258.
10. Barzilai, J., *On the decomposition of value functions11Research supported in part by NSERC.* Operations Research Letters, 1998. **22**(4): p. 159-170.
11. Whitaker, R., *Criticisms of the Analytic Hierarchy Process: Why they often make no sense.* Mathematical and Computer Modelling, 2007. **46**(7-8): p. 948-961.
12. Gass, S.I., *Model world: The great debate - MAUT versus AHP.* Interfaces, 2005. **35**(4): p. 308-312.
13. Smith, J.E. and D.v. Winterfeldt, *Anniversary Article: Decision Analysis in Management Science.* Management Science, 2004. **50**(5): p. 561-574.
14. Forman, E.H., *Facts and fictions about the analytic hierarchy process.* Mathematical and Computer Modelling, 1993. **17**(4–5): p. 19-26.
15. Belton, V. and T. Gear, *On a short-coming of Saaty's method of analytic hierarchies.* Omega, 1983. **11**(3): p. 228-230.
16. Harker, P.T. and L.G. Vargas, *The Theory of Ratio Scale Estimation: Saaty's Analytic Hierarchy Process.* Management Science, 1987. **33**(11): p. 1383-1403.
17. Tu, J. and Z. Wu, *Analytic hierarchy process rank reversals: causes and solutions.* Annals of Operations Research, 2023.
18. Yuen, K.K.F., *Analytic hierarchy prioritization process in the AHP application development: A prioritization operator selection approach.* Applied Soft Computing Journal, 2010. **10**(4): p. 975-989.

19. Chu, A.T.W., R.E. Kalaba, and K. Spingarn, *A comparison of two methods for determining the weights of belonging to fuzzy sets.* Journal of Optimization Theory and Applications, 1979. **27**(4): p. 531-538.
20. Crawford, G. and C. Williams, *A note on the analysis of subjective judgment matrices.* Journal of Mathematical Psychology, 1985. **29**(4): p. 387-405.
21. Lin, C.C., *An enhanced goal programming method for generating priority vectors.* Journal of the Operational Research Society, 2006. **57**(12): p. 1491-1496.
22. Mikhailov, L., *A fuzzy programming method for deriving priorities in the analytic hierarchy process.* Journal of the Operational Research Society, 2000. **51**(3): p. 341-349.
23. Gass, S.I. and T. Rapcsák, *Singular value decomposition in AHP.* European Journal of Operational Research, 2004. **154**(3): p. 573-584.
24. Wang, Y.-M., C. Parkan, and Y. Luo, *Priority estimation in the AHP through maximization of correlation coefficient.* Applied Mathematical Modelling, 2007. **31**(12): p. 2711-2718.
25. Wang, Y.-M., C. Parkan, and Y. Luo, *A linear programming method for generating the most favorable weights from a pairwise comparison matrix.* Computers & Operations Research, 2008. **35**(12): p. 3918-3930.
26. Kou, G. and C. Lin, *A cosine maximization method for the priority vector derivation in AHP.* European Journal of Operational Research, 2014. **235**(1): p. 225-232.
27. Yuen, K.K.F., *The Least Penalty Optimization Prioritization Operators for the Analytic Hierarchy Process: A Revised Case of Medical Decision Problem of Organ Transplantation.* Systems Engineering, 2014. **17**(4): p. 442-461.
28. Zhang, J., et al., *Estimating priorities from relative deviations in pairwise comparison matrices.* Information Sciences, 2021. **552**: p. 310-327.
29. Wang, H., Y. Peng, and G. Kou, *A two-stage ranking method to minimize ordinal violation for pairwise comparisons.* Applied Soft Computing, 2021. **106**: p. 107287.
30. Mangan, P., et al., *Analytic Hierarchy Process Based Land Suitability for Organic Farming in the Arid Region.* Sustainability, 2022. **14**(8): p. 4542.
31. Sekar, K.R., et al., *An improved ranking methodology for malignant carcinoma in multicriterian decision making using hesitant VIKOR fuzzy.* Theoretical Computer Science, 2022. **929**: p. 81-94.
32. Golany, B. and M. Kress, *A multicriteria evaluation of methods for obtaining weights from ratio-scale matrices.* European Journal of Operational Research, 1993. **69**(2): p. 210-220.
33. Mikhailov, L. and M.G. Singh. *Comparison analysis of methods for deriving priorities in the analytic hierarchy process*. in *IEEE SMC'99 Conference Proceedings. 1999 IEEE International Conference on Systems, Man, and Cybernetics (Cat. No.99CH37028)*. 1999.
34. Choo, E.U. and W.C. Wedley, *A common framework for deriving preference values from pairwise comparison matrices.* Computers & Operations Research, 2004. **31**(6): p. 893-908.
35. Gyani, J., A. Ahmed, and M.A. Haq, *MCDM and Various Prioritization Methods in AHP for CSS: A Comprehensive Review.* IEEE Access, 2022. **10**: p. 33492-33511.
36. Yuen, K.K.F., *igm4ahp: Inverse Gram Matrix Methods for Analytic Hierarchy Process*. 2026, Zenodo.https://doi.org/10.5281/zenodo.24223557, https://github.com/kkfyuen/igm4ahp and https://kkfyuen.github.io/igm4ahp.
37. Bozóki, S., *Solution of the least squares method problem of pairwise comparison matrices.* Central European Journal of Operations Research, 2008. **16**(4): p. 345-358.
38. Carrizosa, E. and F. Messine, *An exact global optimization method for deriving weights from pairwise comparison matrices.* Journal of Global Optimization, 2007. **38**(2): p. 237-247.
39. Blankmeyer, E., *Approaches to consistency adjustment.* Journal of Optimization Theory and Applications, 1987. **54**(3): p. 479-488.
40. Johnson, S.G. *The NLopt nonlinear-optimization package*. July 2023; Available from: https://nlopt.readthedocs.io/en/latest/.